\documentclass[12pt]{article}
\usepackage{amsfonts}
\usepackage{amsmath}
\usepackage{epsfig}

\title{A Textbook Case of Pentagram Rigidity}
\author{Richard Evan Schwartz \thanks{Supported by N.S.F. Grant DMS-2102802}}

\newtheorem{theorem}{Theorem}[section]

\newtheorem{lemma}[theorem]{Lemma}

\newtheorem{conjecture}[theorem]{Conjecture}

\def\startproof{{\bf {\medskip}{\noindent}Proof: }}

\def\endproof{$\spadesuit$  \newline}

\def\C{\mbox{\boldmath{$C$}}}%
\def\P{\mbox{\boldmath{$P$}}}%
\def\R{\mbox{\boldmath{$R$}}}%
\def\Z{\mbox{\boldmath{$Z$}}}%

\begin{document}
\maketitle
\begin{abstract}
  In this paper I will explain a rigidity conjecture
  that intertwines the deep diagonal pentagram maps
  and Poncelet polygons. I will also establish a
  simple case of the conjecture, the one involving
  the $3$-diagonal map on a convex $8$-gon with
  $4$-fold rotational symmetry.  This case involves
  a textbook analysis of a pencil of elliptic curves.
  \end{abstract}

\section{Introduction}

\subsection{Conjecture and Result}

Let $\R\P^2$ denote the real projective plane.
A polygon in $\R\P^2$
is {\it convex\/} if its image under a suitable projective
transformation is a convex polygon in the standard affine
patch of $\R\P^2$.  See \S \ref{basics} for definitions.
We call a convex polygon {\it Poncelet\/} if it is inscribed
in one ellipse and circumscribed about another ellipse.
More generally, a {\it Poncelet polygon\/} is one whose
vertices lie in one conic section and whose edges
lie in lines tangent to another conic section.

Let $(n,k)$ be a pair of integers, not both even,
with $n \geq 7$ and $k \in (2,n/2)$.
Given an $n$-gon $P_1$, we let
$P_2=T_k(P_1)$ be the $n$-gon obtained by
intersecting the successive $k$-diagonals of $P_1$.
Figure 1 shows this for $(n,k)=(8,3)$.  
The map $T_k$ is generically defined and
invertible.
The same construction works in 
any field, but convexity is important for us here.
The maps $T_k$ and $T_k^{-1}$ are
always defined on convex $n$-gons, though
the image of a convex $n$-gon under one of these
maps need not be convex.

\begin{center}
\resizebox{!}{1.7in}{\includegraphics{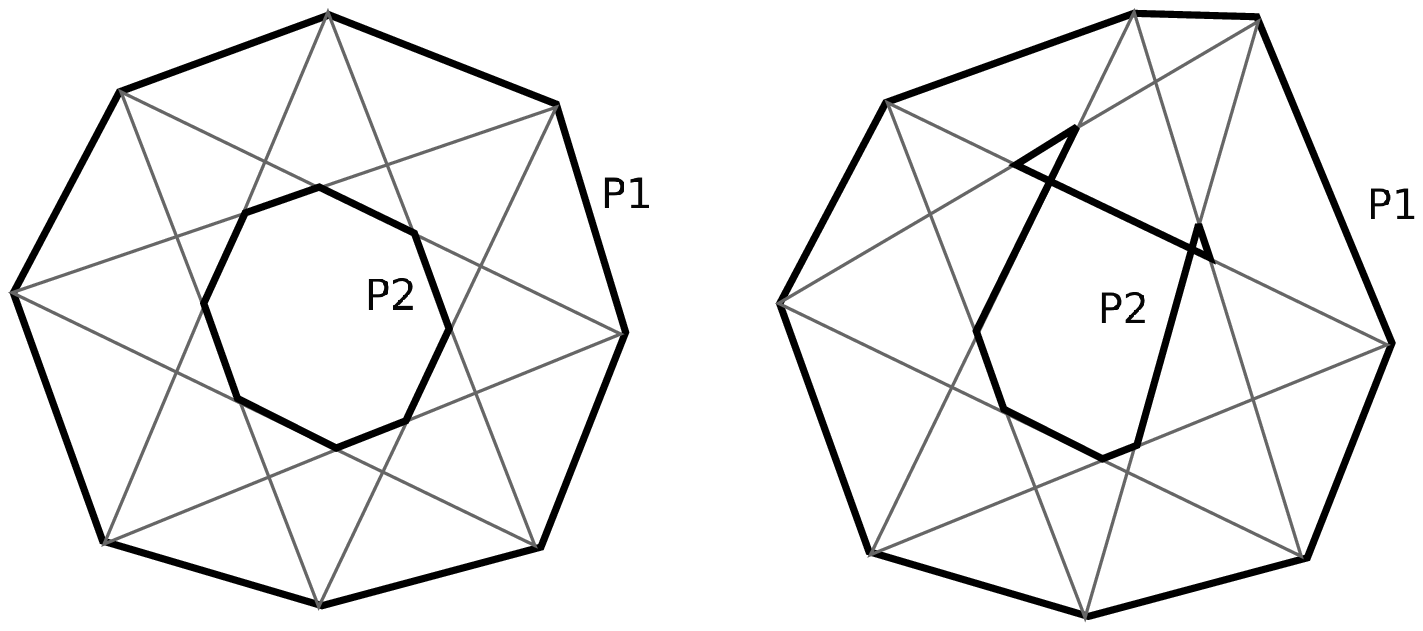}}
\newline
Figure 1: $P_1$ and $P_2=T_3(P_1)$.
\end{center}

Figure 1 gives an example
of where $P_1$ is convex but $P_2$ is not.
Starting with $P_0$ we define the $(n,k)$-{\it pentagram orbit\/}
$\{P_j\}$ where $P_j=T_k^j(P_0)$.

\begin{conjecture}
  Suppose that $\{P_j\}$ is an $(n,k)$-pentagram orbit and
  that $P_j$ is convex for all $j \in \Z$.
  Then $\{P_j\}$ is convex Poncelet for all $j \in \Z$.
\end{conjecture}

I proved in \cite{SCH3} that if $P_0$ is a Poncelet polygon,
$T(P_0)$ and $P_0$ are projectively equivalent.  Thus,
to prove the conjecture it is enough to prove that the
hypotheses force $P_0$ to be convex Poncelet.

In this paper I will prove a simple but nontrivial case of the conjecture.

\begin{theorem}
  \label{octa}
  Suppose that $P_0$ is an $8$-gon with $4$-fold rotational
  symmetry and $\{P_j\}$ is the $(8,3)$-pentagram orbit. Then
  $P_0$ is regular if and only if $P_j$ is convex for all $j$.
  More precisely
  \begin{enumerate}
  \item $P_0$ has $8$-fold dihedral symmetry, with symmetry
    lines containing the vertices, if and only if
    $P_j$ is convex for all $j \geq 0$.
  \item $P_0$ has $8$-fold dihedral symmetry, with symmetry
    lines bisecting the edges, if and only if
    $P_j$ is convex for all $j \leq 0$.
  \end{enumerate}
\end{theorem}

A version of Theorem \ref{octa} appears to be
true when $(8,3)$ is replaced by a general
pair $(n,k)$ when $n$ is even and $k$ is odd and
the $n$-gons have $(n/2)$-fold rotational symmetry.
The proof should be similar.  A much more
interesting generalization is to
the case of centrally symmetric octagons, because these
include all Poncelet octagons. I have just finished
proving this generalization in a much longer paper.
See \cite{SCH4}.

The conjecture is not true for $n$ and $k$ both even.
In this case, $T_k$ is $2$-periodic (modulo scaling) when restricted
to the space of $n$-gons with $(n/2)$-fold rotational
symmetry, and $T_k$ has some non-regular
convex fixed points modulo scale.

\subsection{Context}

The classic case of the pentagram map is $(n,2)$ for
$n \geq 5$.  The case $n=5$ has been
studied e.g. by Clebsch in the $19$th century and
Motzkin \cite{MOT} in middle of the $20$th century.  In $1992$
I wrote a paper \cite{SCH1} defining the pentagram map for
general $n$-gons and proving in the convex case that
the forward orbit shrinks to a point.
Very recently, M. Glick \cite{GLICK1} found
a kind of formula for this collapse point.

It is nice to consider the map $T_2$
as defined on the space ${\cal P\/}_n$ of $n$-gons modulo
projective transformation.  With the correct labeling,
$T_2$ is the identity on ${\cal P\/}_5$ and
$2$-periodic on ${\cal P\/}_6$. I observed
experimentally that the orbits of $T_2$
in general seem to lie on tori.
Motivated (for some reason) by the scattering
transform for the KdV equation, I found \cite{SCH2}
about $n$ algebraically independent invariants
for $T_2$.  These invariants are now called
the {\it monodromy invariants\/} or the
{\it pentagram integrals\/}.

In \cite{OST1}, V. Ovsienko, S. Tabachnikov, and I showed
that $T_2$ has an invariant Poisson structure of
corank $2$ in the odd case and corank $4$ in the even case,
and that the monodromy invariants Poisson-commute with
respect to this structure. This established the complete
Arnold-Liouville integrability of the pentagam map
on the larger space ${\cal T\/}_n$ of so-called
{\it twisted\/} $n$-gons.  Essentially what this means
is that the space ${\cal T\/}_n$, a space of
dimension $2n$, has a singular
foliation by manifolds of dimension about $n$ such that
the restriction of $P_2$ to each manifold is a translation
in suitable coordinates.  When these manifolds are
compact they are necessarily finite unions of tori.

Subsequently, we proved in \cite{OST2}
that $T_2$ is Arnold-Liouville integrable on ${\cal P\/}_n$, which
is naturally a codimension $8$ subvariety of ${\cal T\/}_n$.
For the subset of convex $n$-gons, the manifolds in
the singular foliation are compact and hence finite
unions of tori.
At the same time, F. Soloviev \cite{SOL} proved that the
pentagram map is algebro-geometrically integrable on ${\cal P\/}_n$.
This implies in particular that the torus foliation discussed
above is naturally an abelian fibration, with the
individual tori having natural desciptions as Jacobian
varieties for certain Riemann surfaces.  Very recently,
M. Weinreich \cite{WEIN} proved that $T_2$ is algebro-geometrically
integrable in any field of characteristic not equal to $2$.
In  \cite{GLICK2},  M. Glick related the pentagram
map to a cluster algebra.

By now there are many generalizations of the pentagram
map, and also a number of ways to generate invariant functions
and the invariant Poisson structure.
In \cite{GSTV}, M. Gekhtman, M. Shapiro, S. Tabachnikov,
A. Vainshtein generalized the pentagram map to similar
maps using longer diagonals, and defined on spaces
of so-called {\it corrugated polygons\/} in higher
dimensions.  The work in \cite{GSTV} also generalizes
Glick's cluster algebra and establishes the
complete integrability of these maps in some form.
In \cite{MB1}, G. Mari-Beffa defines higher
dimensional generalizations of the pentagram map
and relates their continuous limits to
various families of integrable PDEs. See also \cite{MB2}.
In \cite{KS}, B. Khesin and F. Soloviev
obtain definitive results about
higher dimensional analogues of the pentagram map, 
their integrability, and their connection to 
KdV-type equations.

The little survey above is not meant to be complete.
Now let me explain how these various results are
related to the Pentagram Rigidity Conjecture above.
First of all, the map $T_k$ is the one used
in \cite{GSTV}.  It would be nice if one could
conclude from \cite{GSTV} that $T_k$ is completely
integrable on ${\cal P\/}_n$, but this has not been
directly worked out.  The spaces
of corrugated polygons are somewhat different
than spaces of ordinary polygons, though in
some sense closely related. Ordinary polygons
are limits of corrugated polygons under a
kind of flattening operation. Let me just
leave it by saying that $T_k$ is certainly
{\it believed\/} to be completely integrable
on ${\cal P\/}_n$ in some sense.

The Pentagram Rigidity Conjecture is really about
the global geometry of the torus foliation
of ${\cal P\/}_n$ (presumably) associated to $T_k$.
The smaller space ${\cal C\/}_n$ of convex $n$-gons
modulo projective transformations is a subset
of ${\cal P\/}_n$. For $k \in [3,n/2)$
  the tori in this foliation probably
are not contained in ${\cal C\/}_n$.  So, if
$T_k$ is not the identity on one of these tori, and
moreover the intersection of the torus
with ${\cal C\/}_n$ is not too large, then
the orbit of $T_k$ on this torus
cannot stay in ${\cal C\/}_n$.  (See
Lemma \ref{minor} below.)  This observation
would prove the conjecture for the
$n$-gons corresponding to this torus.
The Poncelet
polygons in ${\cal C\/}_n$ also lie on these
tori, but $T$ is the identity there.

Motivated by the Pentagram Rigidity Conjecture,
A. Izosimov \cite{IZOS} has recently proved that
if $n$ is odd and $P \in {\cal C\/}_n$ is a fixed point
of $T_2$ then in fact $P$ is a Poncelet polygon.
  So, we now can say that for $n$ odd
points in ${\cal C\/}_n$ are fixed by $T_2$ if and only
if they are Poncelet.  The convexity is important here.
Izosimov gave some easy examples of $n$-gons in
$\C\P^2$ that are fixed by $T_2$ but not Poncelet.
The parity of $n$ is also important.  As
I mentioned above, the result is not true
when $n$ is even.
Presumably, Izosimov's result would also work
for general pairs $(n,k)$ where both numbers are not even.

In the case I consider, that of $T_3$ acting
on $8$-gons with $4$-fold
rotational symmetry, there is just a single
invariant for the map, and its level curves
are nonsingular elliptic curves except when they
contain points corresponding to octagons having
$8$-fold dihedral symmetery.  These are
the curves I analyze in \S \ref{symm}.
These elliptic indeed stretch
outside ${\cal C\/}_8$ in the appropriate sense,
and this is enough to prove Theorem \ref{octa}.
See Figure 2 in \S \ref{proofend}.

The various invariant-generating machines for
the pentagram map probably would turn up
the invariant I found, but these machines
are better developed for
$T_2$ than they are for $T_3$.  I 
just guessed the invariant for $T_3$ by
looking at the picture, and then checked
algebraically that it works.
I will explain in \S \ref{cases} what led me to the invariant.

There are two other connections I want to make
between the Pentagram Rigidity Conjecture and
other areas of mathematics.  When I originally
thought of this conjecture, about $30$ years ago,
I had imagined it as a projective geometry
analogue of the {\it circle packing rigidity\/}
theorem \cite{RS} of B. Rodin and D. Sullivan.
Much
more recently, it occured to me that the
conjecture is something like a discrete analogue
of the Birkhoff-Poritsky Conjecture about
billiards in strictly convex ovals.  This conjecture
says roughly that if a neighborhood of the
boundary of the (cylindrical) billiard phase
space is foliated by invariant curves (corresponding
to caustics) then the oval is an ellipse.

  \subsection{Organization}

  In \S 2 I will give some background
  information about projective geometry
  and also analyze the family of
  elliptic curves that arises in the proof
  of Theorem \ref{octa}.  I will also present
  a few well-known results about complex tori.
  In \S 3 I give the proof of Theorem \ref{octa}.

  The interested reader can download the
  computer program I wrote, which does
  experiments with the $3$-diagonal map
  on centrally symmetric octagons.
  The location of the program is \newline
  {\bf http://www.math.brown.edu/$\sim$res/Java/OCTAGON.tar:\/}

\subsection{Acknowledgements}

I would like to thank Misha Bialy,
Misha Gehktman, Anton Izosimov, Joe Silverman,
Sergei Tabachnikov, and Max Weinreich for
helpful conversations.

\newpage

\section{Preliminaries}

\subsection{Projective Geometry}
\label{basics}

The {\it real projective plane\/}
$\R\P^2$ is the space of lines through
the origin in $\R^3$.  Equivalently it
is the space of scale-equivalence-classes
of nonzero vectors in $\R^3$.
Points in $\R\P^2$ will be denote by
$[x:y:z]$.  This point represents the
line through the origin and $(x,y,z)$.
The quotient map $\R^3 \to \R\P^2$ is called
{\it projectivization\/}.

There is a natural inclusion
$\R^2 \to \R\P^2$ given by
\begin{equation}
  \label{include}
  (x,y) \to [x:y:1].
  \end{equation}
The image of this inclusion is known as the
{\it standard affine patch\/}. I often
identify $\R^2$ with its image under this inclusion, and
when speaking about points in the affine patch
I will often write
$(x,y)$ for $[x:y:1]$.
The inclusion in Equation \ref{include}
has an inverse, given by
\begin{equation}[x:y:z] \to (x/z,y/z).\end{equation}

The {\it line at infinity\/} is the subset of
$\R\P^2$ outside the standard affine patch.
The line at infinity consists of points
of the form $[x:y:0]$.
More generally,
a {\it line\/} in $\R\P^2$ is the set of members
represented by lines in a $2$-dimensional subspace
of $\R^3$.  We can also represent lines by
triples $[a:b:c]$.  This point represents the
linear subspace given by the equation $ax+by+cz=0$.
Conveniently, the line through $2$ points
is represented by the
cross product of the corresponding vectors.
Likewise, the intersection of $2$ lines
is given by the cross product of the corresponding vectors.
These facts make computations with the pentagram map
very easy.

The {\it dual projective space\/}
$\R\P^2_*$ is the space of lines
in $\R\P^2$.   As our notation suggests,
there is an isomorphism between
$\R\P^2$ and $\R\P^2_*$.  The
isomorphism sends the point
represented by $[u:v:w]$ to the
line represented by $[u:v:w]$.
This isomorphism sends collinear points to
coincident lines.

Note that all the same words apply with the
field $\R$ replacing the field $\C$. Thus
$\C\P^2$ is the {\it complex projective plane\/}.
A {\it projective variety\/} in $\C\P^2$ is the
projectivization of the set $V(x,y,z)=0$ where
$V(x,y,z)$ is a homogeneous polynomial in $3$ variables.
This variety is called {\it nonsingular\/} if the (formal) gradient
$\nabla V$ is everywhere nonzero on the set
$V(x,y,z)=0$.  When $V$ is cubic and nonsingular the
corresponding projective variety is a smooth Riemann
surface of genus $1$, also known as a complex torus.
See \cite{SIL}.

\subsection{A Family of Cubics}
\label{symm}

In this section we study the solutions to the
equation
\begin{equation}
  \label{elliptic}
  \frac{(x-y)(x^2+y^2-1)}{xy}=\lambda, \hskip 30 pt
(x-y)(x^2-y^2-1)-\lambda xy =0.
\end{equation}
The second equation is a rearrangement of the first one.
To bring this equation into the form we
mentioned at the end of the last section,
we expand it out and then
homogenize it by padding the $z$-variable.
This gives us the equation

\begin{equation}
V(x,y,z)=x^3 - y^3 - x^2 y + xy^2 - xz^2 + yz^2 - \lambda xyz=0.
\end{equation}

Let $E_{\lambda}$ denote the complex projective variety
corresponding to $V=0$.  Let $\rho$ be reflection in the
line $\{y=-x\}$.  Call this line $L$.  Call a subset
of $\R\P^2$ {\it bounded\/} if it lies in $\R^2$ and
otherwise {\it unbounded\/}.
Below, I will prove two results:
\begin{enumerate}
\item For all $\lambda \not =0, \pm 2, \pm 4 i \sqrt 2$
  the variety $E_{\lambda}$ is nonsingular, and
  hence a complex torus.
\item When $\lambda \in \R-\{-2,0,2\}$, the set
  $E_{\lambda} \cap \R\P^2$ consists of $2$ smooth loops,
  both $\rho$-invariant, one bounded and one unbounded.
  The bounded loop intersects $L$ twice and the unbounded
  loop intersects $L$ once.
\end{enumerate}
Figure 2 in the next chapter shows a rough but topologically accurate picture
of $E_{\lambda} \cap \R^2$ for $\lambda \in (0,2)$.
\newline
\newline
{\bf First Statement:\/}
We want to see that the gradient never vanishes on the level set $V=0$.
We compute
\begin{equation}
\nabla V=
\left[\begin{matrix} V_x \cr V_y \cr V_z \end{matrix} \right]=
\left[\begin{matrix}
-\lambda y z+3 x^2-2 x y+y^2-z^2 \cr -\lambda x z-x^2+2 x y-3 y^2+z^2 \cr 2 z (y-x)-\lambda x y\end{matrix}\right]
\end{equation}

To analyze this, let us first consider the points in the line at infinity that
belong to $V=0$.  Indpendent of $\lambda$, these are the $3$ points
$$[1:1:0], \hskip 30 pt [i:1:0], \hskip 30 pt [-i:1:0].$$
Since both $x,y \not =0$ and $z=0$ for these points, we see
from the equation that the third coordinate of $\nabla V$ is nonzero.
This takes care of these points.

To consider the remaining points of $V$ we can set $z=1$.
If $\nabla V=0$ we have $V_x+V_y=0$. This gives one of
two equations:
\begin{equation}
y=-x, \hskip 30 pt y=\frac{2x-\lambda}{2}.
\end{equation}
When $y=-x$ we have
\begin{equation}
\nabla V = 
\left[\begin{matrix}
\lambda x+6 x^2-1 \cr
-\lambda x-6 x^2+1 \cr
x (\lambda x-4)\end{matrix}
\right]
\end{equation}
This can only vanish when $x=4/\lambda$.
But then
$$V\bigg(\frac{4}{\lambda},-\frac{4}{\lambda},1\bigg)=\frac{256}{\lambda^3} + \frac{8}{\lambda},$$
and this vanishes only if $\lambda = \pm 4 i \sqrt 2$.

When $y=(2x-\lambda)/2$ we have
$$V(x,y,1)=\frac{\lambda(4-\lambda^2)}{8}.$$
This can only vanish when $\lambda=0,\pm 2$.
\endproof

\noindent
{\bf Second Statement:\/}
Let $\lambda \in \R-\{-2,0,2\}$.
The set $E_{\lambda} \cap \R\P^2$ is a finite
disjoint union of smooth loops, permuted by $\rho$.

If $C$ is a bounded component and $\rho(C) \not = C$
then $C \cup \rho(C)$ would intersect some line $4$
times, a contradiction.  Hence $\rho$ preserves each
bounded component, and each bounded component intersects
$L$ twice at right angles.
Since $E_{\lambda}$ intersects $\C\P^2-\C^2$ three times,
and exactly one of these intersection points,
namely $[1:1:0]$,
lies in $\R\P^2$, we see that $E_{\lambda} \cap \R\P^2$
has one unbounded component.

Note that $E_{\lambda} \cap L$ always consists of
$3$ points, namely $(x,-x)$ for $x=0$ and
$$x= =\frac{-\lambda  \pm \sqrt{\lambda^2+32}}{8} \not =0.$$
We conclude that $E_{\lambda} \cap \R\P^2$ must have exactly
two components, both $\rho$-invariant, one
bounded and one unbounded, and the intersections are
as claimed.
\endproof

\subsection{Uniformization}
\label{uniform}

Let $E \subset \C\P^2$ be nonsingular
cubic variety.  As we mentioned above,
$E$ is a complex torus.  Let $f: E \to E$ be
some birational map which is also invertible.
The birational nature of $f$ implie that
all the singularities of $f$ on $E$ are
removable.  This means that $f$ is
a biholomorphic map of $E$ and
orientation preserving.

As is well known, there is also
a biholomorphic map
$\phi: E \to \C/\Lambda$
where $\Lambda \subset \C$ is a lattice.
The map $\phi$ conjugates $f$ to
an isometry of $\C/\Lambda$.
Thus, we may simply equip $E$ with
the coordinates coming from $\phi$
and treat $E$ as a flat torus and
$f: E \to E$ as an orientation preserving isometry.
We call this the {\it flat structure\/} on $E$.

\subsection{Minor Subsets of Tori}

In this section I will prove a general lemma
about flat tori.  The only case required
for the proof of Theorem \ref{octa} is that
of the circle $\R/\Z$,
but the general case might be useful for a
more general version of the conjecture.
The general case rather quickly reduces to
the circle case anyway.

Say that a subset $S \subset \R^n/\Z^n$ is
{\it minor\/} if there is some translation $\phi$ of
$\R^n/\Z^n$ such that
$\phi(S) \subset (0,1/2)^n$.  In general,
say that a subset $S$ of a flat torus is
minor if an affine isomorphism from the
flat torus to $\R^n/\Z^n$ carries $S$ to
a minor subset.

\begin{lemma}
  \label{minor}
  Suppose $p \in S \subset Y$ where $S$ is a minor subset
  of the flat torus $Y$.   Suppose that $f: Y \to Y$ is
  a nontrivial translation.  Then the forward orbit
  $\{f^k(p)|\ k>0\}$ is not contained in $S$.
\end{lemma}

\startproof
By affine symmetry it suffices to prove this when
$Y=\R^n/\Z^n$.  The translation $f$ has a nontrivial
action in at least one coordinate.  Let
$f: \R^n/\Z^n \to \R/\Z$ be the projection onto
this coordinate.  By construction $f(S)$ is minor in
$\R/\Z$ and $f$ covers a nontrivial translation of
$\R/\Z$.  This reduction shows that it suffices to
prove our result for $Y=\R/\Z$.  This is what we do.
If $f(p) \not \in S$ then we are done.
Otherwise $|p-f(p)|<1/2$.  But then the
forward orbit of $p$ is at least $1/4$-dense.
This means that every point of $\R/\Z$ is
within $1/4$ of some point in the forward orbit.
In particular, the point $\zeta$ diametrically opposed
from the midpoint of $S$ has this property.
But then the orbit point that is within
$1/4$ of $\zeta$ is disjoint from $S$.
\endproof

\newpage

\section{Proof of the Result}

\subsection{Formulas}

Let $\cal X$ denote the space of labeled
$8$-gons with $4$-fold rotational symmetry
modulo similarities in the plane.
We normalize so that the $4$-fold symmetry
in question is the map
\begin{equation}
  \rho([x:y:z]) = [-y:x:z].
\end{equation}
This map fixes the origin $(0,0)$ in the affine patch
and preserves the whole affine patch. It is just
rotation by $90$ degrees counterclockwise.

One possibility is that $\rho$ cycles the
vertex labels by $2$ and the other possibility
is that $\rho$ cycles the vertex labels by
$-2$.  We only consider the first possibility;
the second possibility is essentially treated by
symmetry.
For the purpose of getting formulas, we ignore
for now the members of $\cal X$ which have
points on the line at infinity.  We call the
remaining members {\it finite\/}.  In other words,
the finite members lie entirely in the standard
affine patch.
\newline
\newline
\noindent
{\bf The Map:\/}
Every finite member of $\cal X$ has a
canonical representative $P(x,y)$ with
vertices
\begin{equation}
  (1,0), \hskip 4 pt
  (x,y),  \hskip 4 pt
  (0,1), \hskip 4 pt
  (-y,x),  \hskip 4 pt
  (-1,0), \hskip 4 pt
  (-x,-y),  \hskip 4 pt
  (0,-1), \hskip 4 pt
  (y,-x).
\end{equation}
Here $(x,y)$ is really $[x:y:1]$, etc.

Expressed in these coordinates, and with a suitable
labeling scheme, the map $T_3$ is given
by $T_3(x,y)=(x',y')$, where

\begin{equation}
  x'=- A x (B - 2xy), \hskip  30 pt
  y'=+ A y (B + 2xy),
\end{equation}

$$A=\frac{\alpha_{10}+\alpha_{20}}{(1+\alpha_{10})(\alpha_{20}+2\alpha_{30} + \alpha_{40} - \alpha_{11} -2 \alpha_{12}+ \alpha_{22})},
\hskip 20 pt
\alpha_{ij}=x^i y^j + y^i x^j.$$

\begin{equation}
  \label{MAIN}
B = \beta_{10} + 2 \beta_{20} + \beta_{30} + \beta_{12}, \hskip 20 pt
\beta_{ij}=x^iy^j-y^ix^j.
\end{equation}

In other words, $T_3$ sends the polygon $P(x,y)$ to the polygon $P(x',y')$.
I computed this map (and everything else in the paper)
using Mathematica \cite{M}.
\newline
\newline
\noindent
{\bf The Invariant:\/}
Define the function
\begin{equation}
  \Psi(x,y)=\frac{(x-y)(x^2-y^2-1)}{xy}.
\end{equation}
A direct calculation shows that $\Psi \circ T_3=\Psi$.
This is the invariant mentioned in the introduction.
In the next section I will explain where it comes from.
\newline
\newline
{\bf Projective Duality:\/}
Each $8$-gon $P$, defined by its vertices, gives rise to
an $8$-gon $P^*$ in the dual space defined by
the successive lines.  The successive ``vertices'' of
$P^*$ are the successive lines extending the edges of $P$.
Using our isomorphism, we get a second polygon
$(P^*)^{\#}$ in $\cal X$.  The operation $P \to (P^*)^{\#}$
is an involution given algebraically by the map
\begin{equation}
D(x,y)=\bigg(-\frac{y \left(x^2-x\!+\!y^2-y\right)}{x \left(x^2-2 x+y^2+1\right)},\frac{y (x+y-1)}{x \left(x^2-2 x+y^2+1\right)}\bigg)
\end{equation}
Direct calculations show
\begin{equation}
\Psi \circ D= \Psi, \hskip 30 pt
DT_3D^{-1}=T_3^{-1}.
\end{equation}
One can also deduce these equations from
abstract properties of projective duality.
I will leave this to the interested reader.
\newline
\newline
{\bf Symmetries and Factorization:\/}
Define
\begin{equation}
\sigma_1(x,y)=(y,x), \hskip 30 pt \sigma_2(x,y)=(-x,-y).
\end{equation}
A direct calculation shows that
\begin{equation}
\sigma_1 T_3 \sigma_1^{-1}=T_3, \hskip 30 pt
\sigma_2 T_3 \sigma_2^{-1}=T_3^{-1}.
\end{equation}
Geometrically, the map $\sigma_2$ swaps
the regular and star-regular $8$-gons.
Beautifully, a calculation shows that
\begin{equation}
T_3=(D \circ \sigma_2)^2.
\end{equation}
In other words $T_3$ is the square of a simpler map.
Readers familiar with the pentagram map will not be surprised
by this kind of factorization.
The map $D \circ \sigma_2$ satisfies the rule
\begin{equation}
\Psi \circ (D \circ \sigma_2)=-\Psi.
\end{equation}
This equation would probably be the quickest way for the
reader to show, without symbolic manipulation, that
$\Psi$ is an invariant for $T_3$.

\subsection{Special Cases}
\label{cases}

I first noticed that $T_3$ behaved nicely on the sets described
in this section.  I then systematically tested Laurent monomials
in the defining functions for these sets and this led me to $\Psi$.
\newline
\newline
{\bf The Coordinate Axes:\/}
First of all
\begin{equation}
  T_3(x,0)=(-x,0), \hskip 30 pt
  T_3(0,y)=(0,-y),
\end{equation}
So, $T_3$ preserves the coordinate axes and is an involution there.
The corresponding octagons look (to me) like the blades of a circular saw.
\newline
\newline
{\bf The Diagonal Line:\/}
Let $\Delta$ denote the diagonal line $x=y$.
The $8$-gon $P_0=P(x,x)$ has $8$-fold dihedral symmetry, with the lines
of symmetry going through the vertices.  To study the
$(8,3)$-pentagram orbit $\{P_j\}$ we compute
\begin{equation}
  T_3(x,x)=(x',x'), \hskip 30 pt x'=\frac{1+x}{1+2x}.
\end{equation}
The map $T_3$ is given by a projective transformation of $\Delta$.
The fixed points are
\begin{equation}
  p_{\pm}=\pm (1/\sqrt 2,1/\sqrt 2).
\end{equation}

The fixed point $p_+$, which corresponds to the regular $8$-gon,
is attracting.  The fixed point $p_-$, which corresponds to the
star-regular $8$-gon, is repelling.
Thus, every orbit on $\Delta$ aside from $p_-$, is attracted to $p_+$
So, if $P_0$ is not regular then $P_j$ is convex for all $j \geq 0$.
The inverse map $T_3^{-1}$ has
$p_-$ as an attracting fixed point and $p_+$ as a
repelling fixed point.  Hence $P_j$ is non-convex for all
$j$ sufficiently negative.
\newline
\newline
{\bf The Unit Circle:\/}
Let $S^1$ denote the unit circle.
Here $S^1$ corresponds to $8$-gons with
$8$-fold dihedral symmetry in which the
lines of symmetry bisect the sides.  These
$8$-gons are dual to the ones on $\Delta$
and indeed the map $D$ defined above
has the property that $D(\Delta)=S^1$.
Thus, the action of $T_3^{-1}$ on
$S^1$ is conjugate to the action of $T_3$ on
$\Delta$.  In particular, if we start with
$C_0$ convex then $C_j$ is convex for all
$j \leq 0$ but $C_{j}$ is non-convex for all $j$ sufficiently positive.

This is a case that the reader can easily
experiment with.  Just draw a ``stop sign''
and see what the $3$-diagonal map does.
\newline
\newline
    {\bf Other Special Orbits:\/}
    The material here is not needed for the proof of
    Theorem \ref{octa} but it is nice.
Let $L_{\pm 2,1}$ denote the diagonal line $y=x \mp 1$ and
let $L_{\pm 2,2}$ denote the circle of radius $1/\sqrt 2$
centered at $(\mp 1/2,\pm 1/2)$.  The union 
$L_{\pm 2} = L_{\pm 2,1} \cup L_{\pm 2}$ is
the level set $\Psi=\pm 2$.  Each of these
two level sets
is the disjoint union of a
diagonal line and a circle.
The map $T_3^2$ preserves each component and acts
there with order $3$.  For example
\begin{equation}
T_3^2(x,x+1)=(x'',x''+1), \hskip 30 pt
x''=\frac{-1-x}{x}.
\end{equation}

\subsection{Nontriviality}

Let $\lambda \in \R-\{-2,0,2\}$.
We know that $E_{\lambda} \cap \R^2$ contains
one bounded loop and one unbounded loop.
Since $T_3$ preserves $E_{\lambda}$ and $\R^2$ and
$\R\P^2$, we see that $f=T_3^2$ preserves
both components of $E_{\lambda} \cap \R\P^2$.
Let $\Upsilon$ be one of these.

\begin{lemma}
  $f$ cannot be the identity on $\Upsilon$.
\end{lemma}

\startproof
Suppose $f$ is the identity on $\Upsilon$.
Since $f$ is a orientation preserving isometry $E_{\lambda}$ in the flat coordinates,
we see that $f$ must be the identity on $E_{\lambda}$.
We saw in \S \ref{symm} that $E_{\lambda}$ intersects the line
$\{y=-x\}$ in $3$ distinct points. Hence $f(x,-x)=(x,-x)$ for two
distinct nonzero
points $(x_1,-x_1)$ and $(x_2,-x_2)$ in $E_{\lambda}$.
Setting the sum of the coordinates of $f(x,-x)$ equal to $0$, we get
\begin{equation}
  \frac{(4 x^2 (-1 - 2 x^2 + x^4 - 6 x^6 + 32 x^{10}))}{((-1 + x^2 + 4 x^4) (1 - 
  2 x^2 + x^4 + 24 x^6 + 16 x^8))}=0.
  \end{equation}
The only nonzero real roots are
$$x= \pm \frac{1}{2} \sqrt{\frac{1}{2}(1+\sqrt{17})}.$$
But the corresponding points
satisfy $\Psi(x_1,-x_1)=-\Psi(x_2,-x_2) \not =0$ so these
points cannot both lie in $E_{\lambda}$.
This is a contradiction.
\endproof

\noindent
    {\bf Remark:\/} The points $(x_1,-x_1)$ and $(x_2,-x_2)$ constructed in the previous proof
    lie on level sets where $f$ has order $2$.  In particular, $f^2$ is the identity
    on these two level sets.
    \newline

We call an orientation preserving isometry of a metric
circle a {\it translation\/}.

\begin{lemma}
  \label{nontriv}
  $f$ is a nontrivial translation of $\Upsilon$ in the flat coordinates.
\end{lemma}

\startproof
Since $f$ is a nontrivial isometry of $\Upsilon$, we see
that $f$ is either a translation or an orientation
reversing isometry on $\Upsilon$.
In the latter case, $f$ must have order $2$ on $\Upsilon$.
If $f$ is orientation reversing on $\Upsilon$
then $f$ is orientation reversing on nearby
level sets.  This means that $f^2$ is the identity
on all nearby level sets.  But then $f^2$ is the identity
on an open subset of $\R^2$.  Since $f$ is a birational
map, this forces $f^2=T_3^4$ to be the identity. This is false.
\endproof

\subsection{The End of the Proof}
\label{proofend}

Let ${\cal C\/} \subset {\cal X\/}$ denote the
subset of convex $8$-gons.  Let ${\cal C\/}_-$ and ${\cal C\/}_+$ respectively
denote the subsets of $\cal C$ corresponding to $8$-gons with
$8$-fold dihedral symmetry with the lines of symmetry
respectively going through the vertices and bisecting the edges.
Note that ${\cal C\/}_+ \cap {\cal C\/}_-$ is precisely the
point representing the regular $8$-gon.

The analysis of the special sets in \S \ref{cases} reduces us
to considering polygons in
${\cal C\/}-{\cal C\/}_+ - {\cal C\/}_-$.
To finish the proof it suffices to show that
when we have $P_0 \in {\cal C\/}-{\cal C\/}_--{\cal C\/}_+$ there
are indices $i<0<j$ such that
$P_i,P_j \not \in {\cal C\/}$.
Our choice of $P_0$ implies that $\Psi(p) \in \R-\{-2,0,2\}$.
Here $p$ is the point in $\cal C$ representing $P_0$.
But then $p \in E_{\lambda} \cap \R^2$, where
$E_{\lambda}$ is the complex torus discussed in \S \ref{symm}.

Figure 2 shows a picture of the relevant sets.
The shaded semidisk is $\cal C$.
The lightly shaded disk is the unit disk.
The sets ${\cal C\/}_+$ and ${\cal C\/}_-$
are the intersection of ${\cal C\/}$ with the
unit circle and with the diagonal line respectively.
The blue curve is a rough but topologically accurate
sketch of $E_{\lambda} \cap \R\P^2$ when
$\lambda \in (0,2)$.  If $\lambda \in (-2,0)$
the picture would be reflected in the line $\{y=x\}$.

\begin{center}
\resizebox{!}{3.5in}{\includegraphics{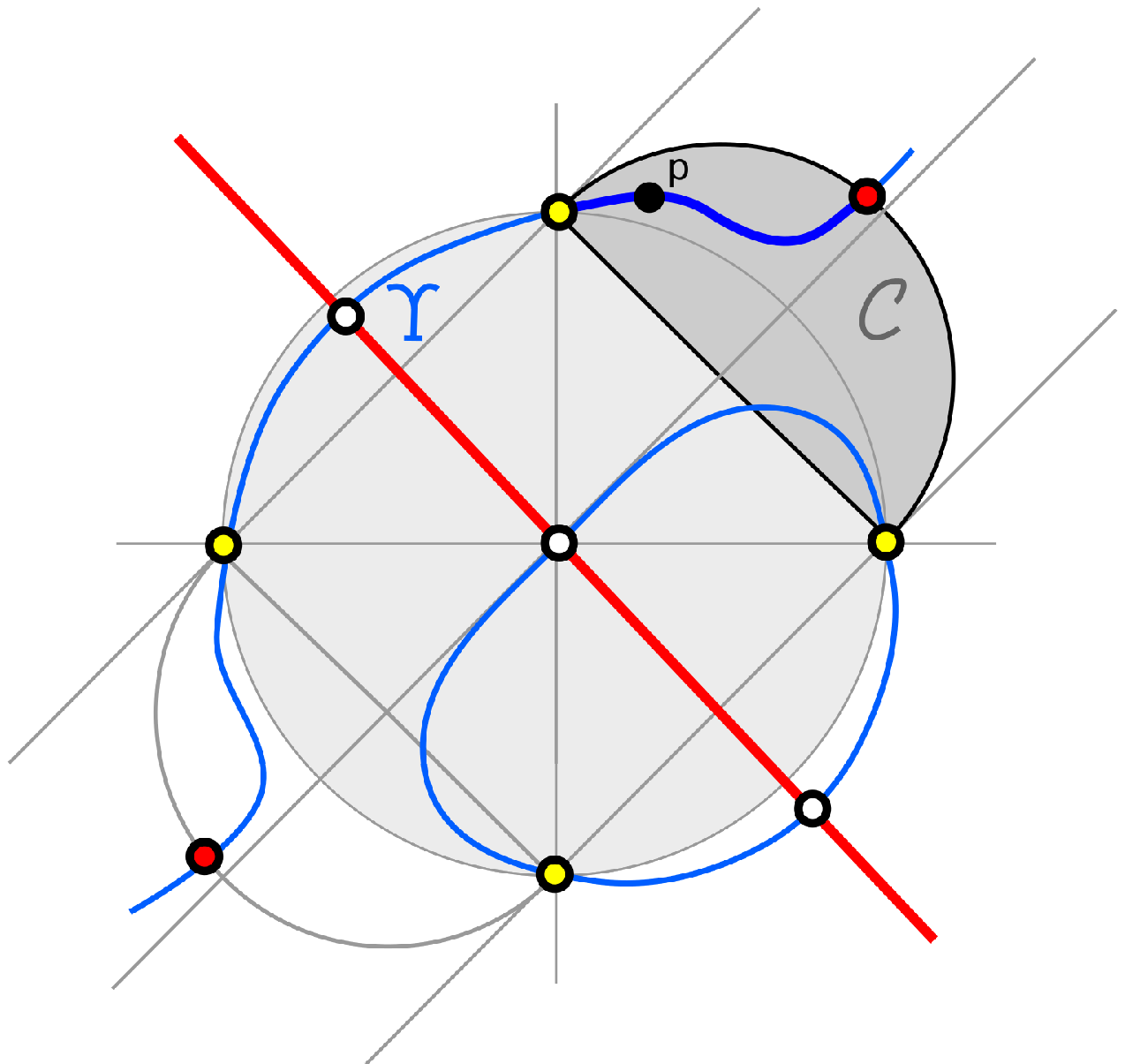}}
\newline
Figure 2: A topologically accurate sketch of $E_{\lambda} \cap \R^2$
for $\lambda \in (0,2)$.
\end{center}

Let $\Upsilon$ be the component of
$E_{\lambda} \cap \R\P^2$ containing $p$.
We equip $\Upsilon$ with its flat coordinates.
Let $f=T_3^2$ as in the previous section.
We know that $f$ is a nontrivial
translation of $\Upsilon$.
Let $\rho$ be the reflection in the
line $\{y=-x\}$.  This is the red diagonal in Figure 2.
We showed in \S \ref{symm} that
$\rho(\Upsilon)=\Upsilon$.  Being holomorphic,
$\rho$ is an isometry of $\Upsilon$ in the flat coordinates.
The two sets
$\Upsilon \cap {\cal C\/}$ and
$\rho(\Upsilon \cap {\cal C\/})$
are disjoint and have equal length.
Hence $\Upsilon \cap {\cal C\/}$ is minor in the
sense of Lemma \ref{minor}.  By Lemma \ref{minor}, we have
$f^j(p) \not \in \Upsilon \cap {\cal C\/}$ for some
$j>0$.  Applying the same argument to $f^{-1}$ gives
some $i<0$ such that $f^i(p) \not \in {\cal C\/}$.
This completes the proof of Theorem \ref{octa}.

\end{document}